\documentclass[10pt]{amsart}
\usepackage{amsmath}
\usepackage{amssymb}
\usepackage{stmaryrd}
\usepackage[all]{xy}
\usepackage{a4wide}
\usepackage[latin1]{inputenc}
\newtheorem{thm}{Theorem}[section]
\newtheorem{cor}[thm]{Corollary}
\newtheorem{lem}[thm]{Lemma}
\newtheorem{prop}[thm]{Proposition}

\newcommand{\To}{\longrightarrow}

\newcommand{\g}{\mathfrak{g}}
\newcommand{\af}{\mathfrak{a}}
\newcommand{\h}{\mathfrak{h}}
\newcommand{\CC}{\mathbb{C}}
\newcommand{\NN}{\mathbb{N}}
\newcommand{\QQ}{\mathbb{Q}}

\newcommand{\ZZ}{\mathbb{Z}}
\newcommand{\indx}[1]{\mathrm{ind}(#1)}
\newcommand{\Ker}[1]{\mathrm{Ker}(#1)}

\newcommand{\Ann}[2]{\mathrm{Ann}_{#1}(#2)}
\newcommand{\ad}[1]{\mathrm{ad}_{#1}}
\newcommand{\SZ}[1]{Z({#1})_s}

\newcommand{\dbl}{\llbracket}
\newcommand{\dbr}{\rrbracket}

    \usepackage{color}
    
\title[Injective hulls of simple modules over Lie superalgebras]
 {Injective hulls of simple modules over finite dimensional nilpotent complex Lie superalgebras} 

\author{Can Hatipo\u{g}lu and Christian Lomp}%
\address{Department of Mathematics, Faculty of Science, University of Porto, Rua Campo Alegre 687, 4169-007 Porto, Portugal}
\email{chatipoglu@alunos.fc.up.pt and clomp@fc.up.pt}
\thanks{Research funded by the European Regional Development Fund through the programme COMPETE and by the Portuguese Government through the FCT - Fundação para a Ciência e a Tecnologia under the project PEst-C/MAT/UI0144/2011. The first author was supported by the grant SFRH/BD/33696/2009}

\begin{document}
\begin{abstract}
We show that the finite dimensional nilpotent complex Lie superalgebras $\g$ whose injective hulls of simple $U(\g)$-modules are locally Artinian are precisely those whose even part $\g_0$ is isomorphic to a
nilpotent Lie algebra with an abelian ideal of codimension $1$ or to a direct product of an abelian Lie algebra and a certain $5$-dimensional or a certain $6$-dimensional nilpotent Lie algebra.
\end{abstract}
\maketitle

\section{Introduction}

Injective modules are the building blocks in the theory of Noetherian rings. Matlis showed that any indecomposable injective module
over a commutative Noetherian ring is isomorphic to the injective hull $E(R/P)$ of some prime ideal $P$ of $R$. He also showed that any injective hull of a simple module is Artinian (see \cite{Matlis} and \cite[Proposition 3]{Matlis_DCC}). 
In connection with the Jacobson Conjecture for Noetherian rings 
Jategaonkar showed in \cite{Jategaonkar_conj} (see also \cite{Cauchon, Schelter}) that the injective hulls of simple modules are locally Artinian provided the ring $R$ is fully bounded Noetherian (FBN). This lead him to answer the Jacobson Conjecture in the afirmative for FBN rings. Recall that a module is called {\it locally Artinian} if every finitely generated submodule of it is Artinian. 
After Jategaonkar's result the question arose whether the condition
\begin{center}$(\diamond)\:\:$ Injective hulls of simple right $A$-modules are locally Artinian\end{center}
was sufficient to prove an affirmative answer of the Jacobson Conjecture which quickly turned out to be not the case. However property $(\diamond)$ remained a subtle condition for Noetherian rings whose meaning is not yet fully understood. Property $(\diamond)$ says that all finitely generated essential extensions of simple right $A$-modules are Artinian. And in case $A$ is right Noetherian property $(\diamond)$ is equivalent to the condition that the class of semi-Artinian right $A$-modules, i.e. modules $M$ that are the union of their socle series, is closed under essential extensions. 

For algebras related to $U(\mathfrak{sl}_2)$ the condition has been examined in \cite{Dahlberg,  injective_modules_over_down-up_algebras, PaulaIan, musson_classification}. One of the first examples of a Noetherian domain that does not satisfy $(\diamond)$ had been found by Ian Musson in \cite{some-examples-of-modules-over-noetherian-rings-musson} concluding that whenever  $\g$ a finite dimensional solvable non-nilpotent Lie algebra, then $U(\g)$ does not satisfy property $(\diamond)$. It is then natural to ask for which finite dimensional complex nilpotent Lie algebras $\g$  its enveloping algebra satisfies $(\diamond)$. We will answer this question completely and will show that those Lie algebras are close to abelian Lie algebras. Slightly more general we can prove our Main Theorem for Lie superalgebras:

\begin{thm}\label{Main_Theorem}
The following statement are equivalent for a finite dimensional nilpotent complex Lie superalgebra  $\g=\g_0\oplus \g_1$:
\begin{enumerate}
\item[(a)] Finitely generated essential extensions of simple $U(\g)$-modules are Artinian.
\item[(b)] Finitely generated essential extensions of simple $U(\g_0)$-modules are Artinian.
\item[(c)] $\indx{\g_0} \geq \dim(\g_0)-2$, where $\indx{\g_0}$ denotes the index of $\g_0$.
\item[(d)] Up to a central abelian direct factor $\g_0$ is isomorphic
\begin{enumerate}
  \item[(i)] to a nilpotent Lie algebra with abelian ideal of codimension $1$;
 \item[(ii)] to the $5$-dimensional Lie algebra $\h_5$ with basis $\{e_1,e_2,e_3,e_4,e_5\}$ and nonzero brackets given by
          $$[e_{1}, e_{2}] = e_{3},\ [e_{1}, e_{3}] = e_{4},\ [e_{2}, e_{3}] = e_{5}.$$
 \item[(iii)] to the $6$-dimensional Lie algebra $\h_6$ with basis $\{e_1,e_2,e_3,e_4,e_5,e_6\}$ and nonzero brackets given by
          $$[e_{1}, e_{2}] = e_{3},\ [e_{1}, e_{4}] = e_{5},\ [e_{2}, e_{4}] = e_{6}.$$
\end{enumerate}
\end{enumerate}
\end{thm}

Together with Musson's  solvable counter example  we have a characterisation of finite dimensional complex solvable Lie algebras $\g$ whose enveloping algebra $U(\g)$ satisfies condition $(\diamond)$.

\begin{cor}
Let $\g$ be a a finite dimensional solvable complex Lie algebra.
$U(\g)$ satisfies $(\diamond)$ if and only if $\g$ is isomporphic up to an abelian direct factor to a Lie algebra with an abelian ideal of codimension $1$ or to $\h_5$ or to $\h_6$.
\end{cor}

The proof of the main Theorem is organized in four steps.
In the first step we show that Noetherian rings whose primitive ideals contain non-zero ideals with a normalizing set of generators satisfy $(\diamond)$, if all of its primitive factors statisfy property $(\diamond)$.
In a second step we verify that ideals of the enveloping algebra $U(\g)$ of a finite dimensional nilpotent Lie superalgebra $\g$ have a  supercentralizing set of generators, which together with the first step shifts our problem to the study of primitive factors of $U(\g)$.
In the third step we combine the description of primitive factors of $U(\g)$ given by  A.Bell and I.Musson as tensor products of the form  $\mathrm{Cliff}_q(\CC)\otimes A_p(\CC)$   with a result of T.Stafford that says that the only Weyl algebra $A_p(\CC)$ satisfying $(\diamond)$ is the first Weyl algebra. A result by E.Herscovich shows that the order $p$ of possible Weyl algebras appearing in the primitive factors of $U(\g)$ is determined by the index $\indx{\g_0}$ of the underlying even part $\g_0$ of $\g$, whch in our case imposes $\indx{\g_0}\geq \dim{\g_0}-2$. The last step lists all finite dimensional nilpotent Lie algebras $\g$ with $\indx{g}\geq \dim(\g)-2$.

\section{Noetherian rings with enough normal elements}

The purpose of this section is to examine the influence that normal elements have on property $(\diamond)$.
Recall that a module $M$ is a subdirect product of a family of modules $\{F_\lambda\}_\Lambda$ if there exists an embedding $\imath: M\rightarrow \prod_{\lambda \in \Lambda} F_\lambda$ into a product of the modules $F_\lambda$ such that for each projection $\pi_\mu:\prod F_\lambda \rightarrow F_\mu$ the composition $\pi_\mu\imath$ is surjective. Compare the next result with \cite[Theorem 1.1]{on-injective-hulls-of-simple-modules-hirano}.

\begin{lem}\label{Diamond_subdirectproduct} A ring $R$ has property $(\diamond)$ if and only if every left $R$-module is a subdirect product of locally Artinian modules.
\end{lem}
\begin{proof} A standard fact in module theory \cite[14.9]{wisbauer} says that every module is a subdirect product of factor modules that are essential extensions of a simple module\footnote{those modules occur in the literature under various names like {\it subdirectly irreducible}, {\it cocyclic}, {\it colocal} or {\it monolithic}}.
Since property $(\diamond)$ is equivalent to subdirectly irreducible modules to be locally Artinian, the Lemma follows.
\end{proof}

A ring extension $R \subseteq S$ is said to be a \textit{finite normalizing extension} if there exists a finite set $\{a_{1}, \ldots a_{k}\}$ of elements of $S$ such that $S = \sum_{i=1}^{k}a_{i}R$ and $a_{i}R = Ra_{i}$, $\forall i=1, \ldots k$. 
The following is an adaption of Hirano's result {\cite[1.8]{on-injective-hulls-of-simple-modules-hirano}}:
\begin{prop}\label{hirano_adapt}
  Let $S$ be a finite normalizing extension of a ring $R$. If $R$ satisfies $(\diamond)$ then so does $S$.
\end{prop}
\begin{proof}
   Let $M$ be a nonzero left $S$-module. By Lemma \ref{Diamond_subdirectproduct} there exists a family $\{N_\lambda \}$ of $R$-submodules of $M$ such that $M/N_\lambda$ is locally Artinian for all $\lambda$ and $\bigcap_\lambda N_\lambda = 0$. For any $R$-submodule $N$ of $M$ denote the largest $S$-submodule of $M$ contained in $N$ by $b(N)$ (called the bound of $N$ in \cite{noncommutative_noetherian_rings}). In fact, $b(N) = \cap_{i=1}^{k}a_{i}^{-1}N$, where $$a_{i}^{-1}N = \{m \in M \mid a_{i}m \in N\}.$$
Since $b(N_\lambda)\subseteq N_\lambda$, we certainly have $\bigcap_{\lambda} b(N_\lambda)=0$.  By \cite[10.1.6]{noncommutative_noetherian_rings}, there is a lattice embedding of $R$-modules $\mathcal{L}(M/b(N_\lambda)) \To \mathcal{L}(M/N_\lambda)$ which implies also that $b(N_\lambda)$ is locally Artinian. Hence $M$ is a subdirect product of locally Artinian $S$-modules.
\end{proof}
As a consequence we have the following.
\begin{cor}\label{tensor_with_fdim_algebra}
 Let $C$ be a finite dimensional algebra and $A$ be any algebra. If $A$ satisfies $(\diamond)$ then $C \otimes A$ satisfies $(\diamond)$ too.
\end{cor}
\begin{proof}
Let $\{x_{1}, \ldots, x_{n}\}$ be a basis of $C$. Then we have $C \otimes A = \sum_{i}^{n}(x_{i} \otimes 1)A$ where each $x_{i} \otimes 1$ is a normal element and so $C \otimes A$ is a finite normalizing extension of $A$ and hence it satisfies $(\diamond)$ by Proposition~\ref{hirano_adapt}.
\end{proof}

A set of elements $\{x_{1}, \ldots, x_{n}\}$ of a ring $R$ is called a \textit{normalizing (resp. centralizing) set} if for each $j = 0,\ldots,n-1$ the image of $x_{j+1}$ in $R/\sum_{i=1}^{j}x_{i}R$ is a normal (resp. central) element.
McConnell showed in \cite{intersection_theorem_for_rings} that every ideal in the enveloping algebra of a finite dimensional nilpotent Lie algebra has a centralizing generator set. In the next section we will show a super version of his result.

\begin{lem}\label{annihilator_of_q}
Let $A$ be a Noetherian algebra, $E$ be a simple $A$-module and $E \leq M$ be an essential extension of left $A$-modules. Let $Q \subseteq\Ann{A}{E}$ be an ideal of $A$ that has a normalizing set of generators. Then $M$ is Artinian if and only if $M' = \Ann{M}{Q}$ is Artinian.
\end{lem}
\begin{proof}
We proceed by induction on the number of elements of the generating set of $Q$. Suppose $Q = \langle x_{1}\rangle$ with $ x_{1} $ being a normal element. Define a map $ f: M \To M $ by $ f(m) = x_{1}m $. This map is $Z(A)$-linear and preserves $A$-submodules of $M$ because if $U \leq M$ is an $A$-submodule of $M$, then $A\cdot f(U) = Ax_{1}U = x_{1}AU = x_{1}U = f(U)$ and so $ f(U) $ if an $A$-submodule of $M$. Since $Q$ is generated by a normal element it satisfies the Artin-Rees property (see \cite[4.1.10]{noncommutative_noetherian_rings}) and so there exists a natural number $n > 0$ such that $ Q^{n}M = x_{1}^{n}M = 0 $. In other words $\Ker{f^{n}} = M$. 
Hence we have a finite filtration
$$0 \subseteq \Ker{f} = \Ann{M}{Q} \subseteq \Ker{f^{2}} \subseteq \ldots \subseteq \Ker{f^{n-1}} \subseteq \Ker{f^{n}} = M$$
whose subfactors are $A/Q$-modules and $f$ induces a submodule preserving chain of embeddings
$$M/\Ker{f^{n-1}} \hookrightarrow \Ker{f^{n-1}} / \Ker{f^{n-2}} \hookrightarrow \ldots \hookrightarrow \Ker{f^{2}} / \Ker{f} \hookrightarrow \Ker{f}.$$
Hence $M$ is Artinian if and only if $M'=\Ker{f}=\Ann{M}{Q}$ is Artinian. Now let $n>0$ and suppose that the assertion holds for all Noetherian algebras and finitely generated essential extensions $E\subseteq M$ of simple left $A$-modules $E$ such that $\Ann{A}{E}$ contains an ideal $Q$ which has a normalizing set of  generators with less than $n$ elements. Let $E\subseteq M$ be a finitely generated essential extension of a simple $A$-module such that $Q\subseteq \Ann{A}{E}$ has a normalizing set of generators $ \{x_{1}, \ldots, x_{n}\} $ of $n$ elements. 
Consider the submodule $M' = \Ann{M}{x_{1}}$. Since $x_{1}$ is a normal element, we can apply the same procedure to conclude that $M$ is Artinian if and ony if $M'$ is Artinian. Let $A' = A / Ax_{1}$ and $ Q' = Q / Ax_{1} $. Then $Q'\subseteq \Ann{A'}{E}$ is generated by the set $ \{\overline{x_{2}}, \ldots , \overline{x_{n}}\} $ of normalizing elements, where $ \overline{x_{i}} $ is the image of $ x_{i} $ in $A'$ for $ i = 2, \ldots,n $. Now, $E \leq M'$ is an essential extension of $A'$-modules such that $Q'E = 0$. Since $Q'$ is generated by a normalizing set of $n-1$ elements, by the induction hypotheses we conclude that $M$ is Artinian if and only if $\Ann{M'}{Q'} = \Ann{M}{Q}$ is Artinian as $A'$-modules and hence also as $A$-modules.
\end{proof}

\begin{lem}\label{Lemma_pre_Conclusion}
Suppose that $A$ is a Noetherian algebra such that every primitive ideal $P$ of $A$ contains an ideal $Q\subseteq P$ which has a normalizing set of generators and  $A/Q$ satisfies $(\diamond)$. Then $A$ satisfies $(\diamond)$.
\end{lem}
\begin{proof}
Let $E$ be a simple $A$-module, $P = \Ann{A}{E}$ and let $E \leq M$ be a finitely generated essential extension of $E$. Let $M' = \Ann{M}{Q}$, where $Q\subseteq P$ is an ideal that has a  normalizing set of generators and with $A/Q$ satisfying $(\diamond)$. Then $E \leq M'$ is a finitely generated essential extension of $A/Q$-modules and so $M'$ is Artinian because $A/Q$ satisfies $ (\diamond) $. Since by Lemma~\ref{annihilator_of_q} $M'$ is Artinian if and only if $M$ is Artinian, it follows that $M$ is Artinian and $A$ satisfies $(\diamond)$.
\end{proof}

Recall that a {\bf superalgebra} is a $\ZZ_2$-graded algebra $A=A_0\oplus A_1$. We denote by $|a|$ the degree of a homogeneous element of $A$. When refering to graded ideals $I$ of $A$ we mean ideals $I=I_0 \oplus I_1$ that are graded with respect to the $\ZZ_2$-grading of $A$. Given any ideal $P$ of $A$ it is easy to see that $Q=P\cap \sigma(P)$ is a graded ideal where $\sigma$ denotes the automorphism :
$$\sigma: A \rightarrow A \qquad a_0+a_1 \mapsto a_0-a_1 \qquad \forall a_0\in A_0, a_1\in A_1.$$

\begin{thm}\label{Conclusion} Let $A$ be a Noetherian superalgebra such that every primitive ideal is maximal and every graded primitive ideal is generated by a normalizing set of generators. Then $A$ satisfies property $(\diamond)$ if and only if every primitive factor of A does.
\end{thm}

\begin{proof}
  The part $(\Rightarrow)$ is clear since the property $(\diamond)$ is inherited by factor rings.\\
  $(\Leftarrow)$ Suppose that every primitive factor of $A$ satisfies $(\diamond)$. Let $E$ be a simple $A$-module, $P = \Ann{A}{E}$, and let $E \leq M$ be an essential extension of $E$. $P$ is maximal by assumption. The ideal $Q=P\cap \sigma(P)$ is graded and has a normalizing set of generators by assumption. If $P$ is graded, then $P=Q$ and $A/Q$ satisfies $(\diamond)$ by hypothesis. 

If $P$ is not graded, then  $P\neq \sigma(P)$ and as $P$ is maximal, $A=P + \sigma(P)$.
Hence $A/Q\simeq A/P \times A/\sigma(P)$. Note that any left $A/Q$-module is $M$ can be written as a direct sum $M=M_1\oplus M_2$ of an $A/P$-module $M_1$ and an $A/\sigma(P)$-module $M_2$. Thus if $E'$ is a simple $A/Q$-module and $M'$ is a finitely generated essential extension of $E'$ as $A/Q$-module, $E'\subseteq M'$ is also a finitely generated essential extension of $A/P$- (resp. $A/\sigma(P)$-) modules. Since $A/P$ satisfies $(\diamond)$ and since $A/P\simeq A/\sigma(P)$, also $A/Q$ satisfies $(\diamond)$.

By  Lemma~\ref{Lemma_pre_Conclusion} we conclude that $A$ satisfies $(\diamond)$. 
\end{proof}

\section{Ideals in enveloping algebras of nilpotent Lie superalgebras}

McConnell showed in \cite{intersection_theorem_for_rings} that every ideal of the enveloping algebra of a finite dimensional nilpotent Lie algebra has a centralizing set of generators. We intend to prove an analogous result for superalgebras. 

The {\bf supercommutator} of two homogeneous elements $a,b$ is the element
$$ \dbl a,b \dbr \ := \ ab - (-1)^{|a||b|}ba$$
and is extended bilinearly to a form $\dbl - , - \dbr : A \rightarrow A$.
The {\bf supercenter} of $A$ is the set $\SZ{A} = \{a\in A \mid \forall b\in A: \dbl a,b \dbr = 0\}$ and its elements are called {\bf supercentral}. Given a supercentral element $a\in A$, the ideal $I=Aa$ is a graded and $A/Aa$ is again a superalgebra. We say that a set of elements $\{x_{1}, \ldots, x_{n}\}$ of a superalgebra $A$ is a \textit{supercentralizing set} if for each $j = 0,\ldots,n-1$ the image of $x_{j+1}$ in $A/\sum_{i=1}^{j}x_{i}A$ is a supercentral element.

A homogenous superderivation of a superalgebra $A$ is a linear map $f: A \To A$ such that
$$f(ab) = f(a)b + (-1)^{|a||b|}af(b)$$
for all homogeneous $a,b \in A$. The supercommutator $\dbl a, - \dbr$ for a homogeneous element is a superderivation. In case $|a|=0$, $\dbl a, - \dbr$ is a derivation of $A$.

Let $\g=\g_0 \oplus \g_1$ be a Lie superalgebra and choose a basis $\{x_1,\ldots, x_n\}$ of $\g_0$ and a basis $\{y_1,\ldots, y_m\}$ of $\g_1$.
The PBW theorem for Lie superalgebras (see \cite{Behr}) says that the monomials $x_1^{\alpha_1}\cdots x_n^{\alpha_n} y_1^{\beta_1}\cdots y_m^{\beta_m}$ with $\alpha_i, \beta_j \in \NN_{0}$ and 
$\beta_i\leq 1$ form a basis of the enveloping algebra $A=U(\g)$. For $i\in \{0,1\}$ let 
$$A_i=\mathrm{span} \{ x_1^{\alpha_1}\cdots x_n^{\alpha_n} y_1^{\beta_1}\cdots y_m^{\beta_m} \mid 
\beta_1+\cdots + \beta_m = i (\mathrm{mod} 2) \}.$$
Then $A=A_0\oplus A_1$ is a superalgebra such that the degree of a homogeneous element of $\g$ equals its degree in $A$.
For any $x\in \g$, the adjoint action of $x$ on $A$ is defined by
$$\ad{x}: A \rightarrow A \qquad \ad{x}(a)=\dbl x,a \dbr \:\:\forall a\in A.$$
By definition of the enveloping algebra we have for all $x,y \in \g$:
$$ \ad{x}(y) = \dbl x,y \dbr \ =\  [x,y].$$

The following Lemma follows from a direct computation which we carry out for the convenience of the reader.
\begin{lem}\label{Lemma1}
For any $x,y \in \g$ one has \begin{equation}\ad{x}\circ \ad{y} - (-1)^{|x||y|} \ad{y}\circ\ad{x} = \ad{[x,y]}.\end{equation}
\end{lem}

\begin{proof} Let $a$ be a homogeneous element of $A$, $x,y \in \g$.
\begin{eqnarray*}	
\lefteqn{\dbl x, \dbl y,a\dbr\dbr-(-1)^{|x||y|}\dbl y, \dbl x,a\dbr\dbr}\\
&=& x(ya-(-1)^{|y||a|}ay)-(-1)^{|x|(|y|+|a|)}(ya-(-1)^{|y||a|}ay)x \\
&& -(-1)^{|x||y|}\{y(xa-(-1)^{|x||a|}ax)-(-1)^{|y|(|x|+|a|)}(xa-(-1)^{|x||a|}ax)y \}\\
&=& xya+(-1)^{|x||y|+|x||a|+|y||a|}ayx -(-1)^{|x||y|}yxa-(-1)^{|a||y|+|x||a|}axy\\
&=& [x,y]a+(-1)^{|a|(|x|+|y|)}a[x,y] = \dbl [x,y], a\dbr
\end{eqnarray*}
\end{proof}

\begin{cor}\label{corollary_lemma1}
For any $x\in \g_1$, then $\ad{x}^2 = 0$.
\end{cor}
\begin{proof}
 Set $y=x$ in Lemma \ref{Lemma1}, then $2\ad{x}^2 = \ad{[x,x]} = 0$.
\end{proof}

Recall that a map $f:A\rightarrow A$ is called locally nilpotent if for every $a\in A$ there exists a number $n(a)\geq 0$ such that $f^{n(a)}(a)=0$.

\begin{prop}\label{all_inner_locally_nilpotent}
 Let $\g$ be a finite dimensional nilpotent Lie superalgebra. Then $\ad{x}$ is locally nilpotent superderivation of $A=U(\g)$, for every homogeneous element $x\in \g$.
\end{prop}

\begin{proof}
 In case $x \in \g_1$ is odd, we see from Corollary \ref{corollary_lemma1} that $\ad{x}$ is nilpotent.
In case $x\in \g_0$ is even, then $\ad{x}=\dbl x, - \dbr$ is an ordinary derivation of $A$. Let $r$ be the nilpotency degree of $\g$, i.e. $\g^{r}=0$.
Then for any $a\in \g$ we have $\ad{x}^r(a)=0$. Let $m\geq 0$. Suppose that for every monomial $a\in A$ of length $m$ there exists $n(a)\geq 0$ such that $\ad{x}^{n(a)}(a)=0$. Let $y\in \g$. Then 

$$\ad{x}^{n(a)+r}(ay) = \sum_{i=1}^{n(a)+r} { {n(a)+r} \choose {i} } \ad{x}^i(a)\ad{x}^{n(a)+r-i}(y) = 0.$$
By induction $\ad{x}$ is locally nilpotent on all basis elements of $A$.
\end{proof}

Given an $l$-tuple of superderivations $\partial=(\partial_1, \ldots, \partial_l)$ of a superalgebra $A$ we say that a subset $X$ of $A$ is $\partial$-stable if $\partial_i(X)\subseteq X$ for all $1\leq i \leq l$. Note that if all superderivations $\partial_i$ are inner, i.e. $\partial_i=\dbl x_i,-\dbr$ for some homogeneous $x_i\in A$, then any ideal $I$ is $\partial$-stable. 

\begin{thm}\label{superalgebra_case}
 Let $A$ be a superalgebra with locally nilpotent superderivations $\partial_1, \ldots, \partial_l$ such that $\bigcap_{i=1}^l \ker\partial_i \subseteq \SZ{A}$ and for all $i\leq j$ there exist $\lambda_{i,j}\in \CC$.
\begin{equation}\label{relation_eq} \partial_i \circ \partial_j - \lambda_{i,j} \partial_j \circ \partial_i \in \sum_{s=1}^{i-1} \CC \partial_s.\end{equation}
Then any non-zero $\partial$-stable ideal $I$ of $A$ contains a non-zero supercentral element. In particular if $I$ is graded and Noetherian, then it contains  a supercentralizing set of generators.
\end{thm}

\begin{proof}
For each $1\leq t \leq l$ set $K_t = \bigcap_{i=1}^t \ker\partial_i$. 
We will first show that $K_i$ are $\partial$-stable subalgebras of $A$.
Let $1\leq t, j \leq l$ and $a\in K_t$. If $j\leq t$, then $\partial_j(a)=0\in K_t$ by definition. Hence suppose $j>t$. By hypothesis for any $1\leq i \leq t<j$ we have
$$\partial_i(\partial_j(a)) = \lambda_{i,j}\partial_j(\partial_i(a)) + \sum_{s=1}^{i-1} \mu_{i,j,s} \partial_s(a) = 0$$
for some $\lambda_{i,j}, \mu_{i,j,s} \in \CC$. Thus $\partial_j(a)\in K_t$.

To show that $I$ contains a non-zero element of the supercentre of $A$ note that since $\partial_1$ is locally nilpotent, for any $0\neq a \in I$ there exists $n_1\geq 0$ such that $0\neq a'=\partial_1^{n_1}(a) \in \ker \partial_1=K_1$. Since $I$ is $\partial_1$-stable, $a'\in I\cap K_1$. Suppose $1\leq t \leq l$ and $0\neq a_t \in I\cap K_t$, then since $\partial_{t+1}$ is locally nilpotent, there exists $n_{t+1}\geq 0$ such that $0\neq a'=\partial_{t+1}^{n_{t+1}}(a) \in \ker \partial_{t+1}$. Since $I$ and $K_t$ are $\partial$-stable, we have $a'\in I\cap K_{t+1}$. Hence for $t=l$, we get $0\neq I\cap K_t \subseteq I\cap \SZ{A}$.

Assume that $I$ is graded and Noetherian and let $0\neq a=a_0+a_1 \in I \cap \SZ{A}$. Since $I$ and $\SZ{A}$ are graded, both parts $a_0$ and $a_1$ belong to $I\cap \SZ{A}$, one of them being non-zero. Thus we might choose $a$ to be homogeneous. Let $J_1=Aa$ be the graded ideal generated by $a$, then all superderivations $\partial_i$ lift to superderivations of $A/J_1$ satisfying the same relation (\ref{relation_eq}) as before. Moreover $I/J_1$ is a graded Noetherian $\partial$-stable ideal of $A/J_1$. Applying the procedure of obtaining a supercentral element to $I/J_1$ in $A/J_1$ yields a supercentral homogeneous element $a'\in I/J_1 \cap \SZ{(A/J_1}$. Set $J_2=Aa + Aa'$. Continuing in this way leads to an ascending chain of ideals $J_1\subseteq J_2 \subseteq \cdots \subseteq I$ that eventually has to stop, i.e. $I=J_m$ for some $m$. By construction, the generators used to build up $J_1, J_2, \ldots, J_m$ form a supercentralizing set of generators for $I$.
\end{proof}

In order to apply the last Proposition to the enveloping algebra of a finite dimensional nilpotent Lie superalgebra $\g$, we have to choose an appropriate basis of homogeneous elements.  Without loss of generality we might assume that $\g$ has a refined central series
$$ \g=\g^n \supset \g^{n-1} \supset \g^{n-2} \supset \cdots \supset \g^1 \supset \g^{0}=\{0\}.$$
with $[\g,\g^i] \subseteq \g^{i-1}$ and $dim(\g^i/\g^{i-1})=1$ for all $1\leq i \leq n$.
Let $x_1,x_2, \ldots, x_{n}$ be a basis of $\g$ such that each element $x_i+\g^{i-1}$ is non-zero (and hence forms a basis) in $\g^i/\g^{i-1}$. Actually each $x_i$ is homogeneous, since if $x_i={x_i}_0 + {x_i}_1$ with ${x_i}_j$ homogeneous, then as ${x_i}_0$ and ${x_i}_1$ cannot be linearly independent as $g^i/g^{i-1}$ is 1-dimensional, one of them belongs to $\g^{i-1}$.

\begin{cor}\label{supercentral}
Any graded ideal of the enveloping algebra of a finite dimensional nilpotent Lie superalgebra has a supercentralizing set of generators.
\end{cor}

\begin{proof}
  Let $\g$ and $A=U(\g)$ be as above, as well as the chosen basis of $\g$ $x_1, \ldots, x_n$ of homogeneous elements. Set $\partial_i=\ad{x_i}$. By Proposition \ref{all_inner_locally_nilpotent} all superderivations $\partial_i$ are locally nilpotent.
Let $i<j$, then $[x_i,x_j] \in \g^{i-1}$ show that there are scalars $\mu_{i,j,s}\in \CC$ such that 
$$[x_i,x_j] = \sum_{s=1}^{i-1} \mu_{i,j,s} x_s.$$
Note that $\ad{[x_i,x_j]} = \sum_{s=1}^{i-1} \mu_{i,j,s} \ad{x_s}.$
Therefore, using Lemma \ref{Lemma1}, we have 
$$\partial_i\circ \partial_j = (-1)^{|x_i||x_j|}\partial_j\circ \partial_i + \sum_{s=1}^{i-1} \mu_{i,j,s} \partial_s.$$
Hence the assumptions of Theorem \ref{superalgebra_case} are fulfilled and our claim follows (since $A$ is Noetherian).
\end{proof}

This last result with Theorem~\ref{Conclusion} gives the following:

\begin{cor}\label{conclusion_super_lie}
Let $\g$ be a finite dimensional nilpotent Lie superalgebra. Then $U=U(\g)$ satisfies property $(\diamond)$ if and only if every primitive factor of $U$ does.
\end{cor}

\begin{proof}
By Corollary \ref{supercentral} any graded ideal is generated by supercentral hence normal elements. Moreover every primitive ideal of $U(\g)$ is maximal by \cite[Corollary 1.6]{Letzter}. Hence the result follows from Theorem \ref{Conclusion}.
\end{proof}

\section{Primitive factors of nilpotent Lie superalgebras}

It is a standard fact that primitive factors of enveloping algebras of  finite dimensional nilpotent Lie algebras are Weyl algebras.
Recall that the $n$th Weyl algebra over $\CC$ is the algebra $A_{n}(\CC)$ generated by $2n$ elements $x_{1}, \ldots , x_{n}, y_{1}, \ldots , y_{n}$ subject to the relations $x_{i}y_{j}-y_jx_i= \delta_{ij}$, for all $1\leq i,j\leq n$.

A.Bell and I.Musson showed in \cite{primitive_factors_of_enveloping_algebras_of_nilpotent_lie_superalgebras-musson_bell} that primitive factors of enveloping algebras of  finite dimensional nilpotent Lie superalgebras are of the form $\mathrm{Cliff}_q(\CC)\otimes A_p(\CC)$ where $\mathrm{Cliff}_q(\CC)$ is a Clifford algebra. 
We know from \cite{Lam_QuadraticForms} that 
$$\mathrm{Cliff}_0(\CC) = \CC,\qquad \mathrm{Cliff}_1(\CC) = \CC \times \CC, \qquad \mathrm{Cliff}_2(\CC) = M_2(\CC)$$
and $\mathrm{Cliff}_{n+2}(\CC) = \mathrm{Cliff}_n(\CC) \otimes M_2(\CC)$ for all $n>2$. The next Lemma shows that property $(\diamond)$ is stable under tensoring with a Clifford algebra:

\begin{lem}\label{Clifford}
A $\CC$-algebra $A$ satisfies $(\diamond)$ if and only if $\mathrm{Cliff}_q(\CC) \otimes A$ satisfies $(\diamond)$ for all (for one) $q$.
\end{lem}

\begin{proof}
By Corollary~\ref{tensor_with_fdim_algebra}, $\mathrm{Cliff}_q(\CC) \otimes A$ satisfies $(\diamond)$ if $A$ does. On the other hand suppose that there exists $q>0$ such that $\mathrm{Cliff}_q(\CC) \otimes A$ satisfies $(\diamond)$. 
If $q=2m$ is even, then $\mathrm{Cliff}_q(\CC) \otimes A = M_{2^m}(A)$ which is Morita equivalent to $A$. Since $(\diamond)$ is a Morita-invariant property as the equivalence between module categories yields lattice isomorphisms of the lattice of submodules of modules, we get that $A$ satisfies $(\diamond)$.
If $q=2m+1$ is odd, then $\mathrm{Cliff}_q(\CC) \otimes A  = M_{2^m}(A)\times M_{2^m}(A)$. Since $A$ is Morita equivalent to the factor $M_{2^m}(A)$ it also satisfies $(\diamond)$. \end{proof}

The question is hence which Weyl algebras do satisfy $(\diamond)$. Being a semiprime Noetherian ring of Krull dimension 1, the first Weyl algebra $A_{1}(\CC)$ satisfies the property $(\diamond)$ \cite{injective_modules_over_down-up_algebras}. However, for $n \geq 2$, the Weyl algebra $A_n=A_{n}(\CC)$ does not satisfy the property $(\diamond)$. In\cite{nonholonomic_modules_over_weyl_algebras_and_enveloping_algebras} J. T. Stafford constructs a simple $A_{n}(\CC)$-module which has an essential extension of Krull dimension $n-1$:
\begin{thm}[T.Stafford {\cite[Theorem 1.1, Corollary 1.4]{nonholonomic_modules_over_weyl_algebras_and_enveloping_algebras}}]\label{stafford's_theorem}
  For $2 \leq i \leq n$ pick $\lambda_{i} \in \CC$ that are linearly independent over $\QQ$. Then the element
  $$\alpha = x_{1} + y_{1}\left(\sum_{2}^{n}\lambda_{i}x_{i}y_{i}\right) + \sum_{2}^{n}(x_{i} + y_{i})$$
generates a maximal right ideal of $A_n=A_{n}(\CC)$. In particular
$A_{n} / x_{1}\alpha A_{n}$ is an essential extension of the simple $A_{n}$-module $A_{n} / \alpha A_{n}$ by the module $A_{n} / x_{1}A_{n}$, which has Krull dimension $n-1$.
\end{thm}
Since Artinian modules are exactly the ones with Krull dimension zero, this implies that $A_{n}(\CC)$ satisfies the property $(\diamond)$ if and only if $n=1$. Stafford's result is a key ingredients in the proof of our main theorem.
The order of Weyl algebras appearing in the primitive factors of enveloping algebras $U(\g)$ of finite dimensional nilpotent Lie superalgebras $\g$ has been determined by E.Herscovich in \cite{dixmier_map_for_nilpotent_super_lie_algebras-herscovich} and is related to the index of the underlying even part of $\g$.

Let $f\in \g^*$ be a linear functional on a Lie algebra $\g$ and set
$$g^{f} = \{ x \in \g \mid f([x,y]) = 0, \: \forall y\in \g \}$$ be the orthogonal subspace of $\g$ with respect to the bilinear form $f([-,-])$.
The number  $$\indx{\g}:=\inf_{f \in \g^{*}} \dim \g^{f}$$ is called the \textit{index} of $\g$.

\begin{thm}[E.Herscovich {\cite{dixmier_map_for_nilpotent_super_lie_algebras-herscovich}}, A.Bell \& I.Musson {\cite{primitive_factors_of_enveloping_algebras_of_nilpotent_lie_superalgebras-musson_bell}}]\label{Estanislao}
Let $\g$ be a finite dimensional nilpotent complex Lie superalgebra.
\begin{enumerate}
\item For  $f\in \g^*$ there exists a graded primitive  ideal $I(f)$ of $U(\g)$ such that
$$U(\g)/I(f)\simeq \mathrm{Cliff}_q(\CC) \otimes A_p(\CC),$$ where
$2p=\dim(\g_0/\g_0^f) \leq \dim(\g_0)-\indx{\g_0}$ and $q\geq 0$.
\item For every graded primitive ideal $P$ of $U(\g)$ there exists $f\in \g^*$ such that $P=I(f)$.
\end{enumerate}
\end{thm}

Combining Stafford's and Herscovich's results with \ref{conclusion_super_lie} leads now easily to the following

\begin{prop}\label{diamond-index}
Let $\g=\g_0\oplus \g_1$ be a finite dimensional nilpotent complex Lie superalgebra. Then $U(\g)$ satisfies $(\diamond)$ if and only if $\indx{\g_0} \geq \dim(\g_0)-2$.
\end{prop}

\begin{proof}
$(\Rightarrow)$ By Theorem~\ref{Estanislao} each primitive factor of $U(\g)$ is of the form $ \mathrm{Cliff}_q(\CC) \otimes A_p(\CC) $ where $2p=\dim(\g_0/\g_0^f) = \dim(\g_0)-\dim{\g_0^f}$. Since the property $(\diamond)$ is inherited by factor rings this implies together with Theorem~\ref{stafford's_theorem} and Lemma~\ref{Clifford} that $p \leq 1$, that is $\dim{g_0^f} \geq \dim(\g_0)-2$ ,i.e. $\indx{\g_0}\geq \dim(\g_0)-2$.\\
$(\Leftarrow)$ If $\indx{\g_0} \geq \dim(\g_0)-2$ then the primitive factors of $U(\g)$ are either of the form $\mathrm{Cliff}_q(\CC)$ or $\mathrm{Cliff}_q(\CC) \otimes A_1(\CC)$. Thus the primitive factors of $U(\g)$ satisfy the property $(\diamond)$ by Lemma~\ref{Clifford}. This implies together with Corollary~\ref{conclusion_super_lie} that $U(\g)$ satisfies $(\diamond)$.
\end{proof}

\section{Nilpotent Lie algebras with almost maximal index}
In this last section we will classify all finite dimensional complex Lie algebras $\g$ with index greater or equal to $\dim{\g}-2$. It is clear that if $\indx{\g}=\dim{\g}$, then $\g$ is abelian.  We say that a Lie algebra $\g$ has {\it almost maximal index} if $\indx{\g}=\dim(\g) - 2$.

As a first step we show that a direct product $\g_1 \times \g_2$ of two Lie algebras $\g_{1}$ and $\g_{2}$ has almost maximal index if and only if one of them is abelian and the other one has almost maximal index. Recall that the Lie bracket of the direct product 
$\g = \g_{1} \times \g_{2}$ is defined as 
$$[(x_{1}, y_{1}), (x_{2}, y_{2})] := ([x_{1}, x_{2}], [y_{1}, y_{2}])$$
for all $x_{1}, x_{2} \in \g_{1},\ y_{1}, y_{2} \in \g_{2}$. For the product algebra, we have the following formula:

\begin{lem}\label{index_formula}
 For Lie algebras $\g_1, \g_2$ the following formula holds: $$\indx{\g_1 \times \g_2} = \indx{\g_1} + \indx{\g_2}.$$
In particular $\g_1\times \g_2$ has almost maximal index if and only if one of the factors has almost maximal index and the other factor is Abelian.
\end{lem}
\begin{proof} Set $\g=\g_1 \times \g_2$. Since $\g^* = \g_1^* \times \g_2^*$, for all $f \in \g^{*}$, we have 
$\dim{\g^f} = \dim{\g_1^{f_1}} + \dim{\g_2^{f_2}}$, with $f_i = f\epsilon_i \in \g_i^*$ and inclusions $\epsilon_i:\g_i\rightarrow \g$.
Thus $\indx{\g} = \indx{\g_1} + \indx{\g_2}$. Note that in general $\indx{\g_i}=\dim(\g_i) - 2n_i$ for some $n_i\geq 0$ and let $\g=\g_1\times \g_2$.
Hence $$\indx{\g}=\indx{g_1}+\indx{g_2} = \dim(\g_1)-2n_1+\dim(\g_2)-2n_2 = \dim(\g) - 2(n_1+n_2) = \dim(\g)-2$$ if and only if $n_1+n_2=1$ which shows our claim.
\end{proof}
The Lemma together with Proposition \ref{diamond-index} implies:

\begin{prop}
  Let $\g$ be a finite dimensional complex nilpotent Lie algebra. 
Then \\$U(\g)[x_{1}, \ldots, x_{n}]$ has the property $(\diamond)$ if and only if $U(\g)$ has the property $(\diamond)$.
\end{prop}
\begin{proof}
  Suppose that $U(\g)$ has the property $(\diamond)$. We have
$$
U(\g)[x_{1}, \ldots, x_{n}]  = U(\g) \otimes \CC[x_{1}, \ldots, x_{n}]
                             = U(\g) \otimes U(\mathfrak{a})
                             = U(\g \oplus \mathfrak{a})
$$
for an $n$-dimensional Abelian Lie algebra $\mathfrak{a}$. Since $U(\g)$ satisfies $(\diamond)$, $\g$ has index at least $\dim(\g) - 2$. By Lemma~\ref{index_formula}, we have $\operatorname{ind}(\g \oplus \mathfrak{a}) \geq \dim(\g) + n - 2 = \dim(\g \oplus \mathfrak{a}) - 2$. Since $\g \oplus \mathfrak{a}$ is nilpotent, it follows by Proposition~\ref{diamond-index} that $U(\g \oplus \mathfrak{a})$ satisfies $(\diamond)$. Thus $U(\g)[x_{1}, \ldots, x_{n}]$ also satisfies $(\diamond)$. Conversely, if the polynomial algebra $U(\g)[x_{1}, \ldots, x_{n}]$ has the property $(\diamond)$, then $U(\g)$ also has it since it is inherited by factor rings.
\end{proof}

Note that in general it seems unknown whether property $(\diamond)$ is inherited by forming polynomial rings. Lemma \ref{index_formula} also shows that we can ignore abelian direct factors for the characterization of Lie algebras with almost maximal index. The following Proposition will classify those Lie algebras.
\begin{prop}\label{listOfDiamond} A finite dimensional nilpotent Lie algebra $\g$ has almost maximal index if and only if $\g$ has an abelian ideal of codimension $1$ or if $\g$ is isomorphic (up to an abelian direct factor) to $\h_5$ or $\h_6$.
\end{prop}

\begin{proof}
Let $\g$ be a finite dimensional nilpotent Lie algebra of dimension $n$ and index $n-2$ and suppose that $\g$ does not have an abelian ideal of codimension $1$.  There exists a functional $f\in \g^*$ such that $\dim\g^{f} = \indx{\g}=n-2$. 
By \cite[1.11.7]{enveloping_algebras},  $\g^{f}$ is an abelian Lie subalgebra of $\g$. By \cite[5.1]{BurdeCeballos} there exists an abelian ideal $\af$ of $\g$ of codimension $2$.
Let $\{e_{1}, \ldots , e_{n}\}$ be a basis of $\g$, such that $\{e_3, \ldots, e_n\}$ is a basis for $\af$.
Since $\af$ is abelian, the matrix of brackets $[e_i,e_j]$ has the form
$$M=([e_i,e_j])_{i,j} = \left(
    \begin{array}{cc}
      \phantom{-}A & B \\
      -B^{t} & \mathbf{0} \\
    \end{array}
  \right)
$$
where $A$ is a $2\times 2$ skew-symmetric matrix, $B$ is a $2\times(n-2)$ matrix with entries in $\af$, and $\mathbf{0}$ is the $(n-2)\times(n-2)$ zero matrix. 
Let $$B_{ij}=\left(\begin{array}{ll} [e_1,e_i] & [e_1,e_j] \\ {[e_2,e_i]} & [e_2,e_j] \end{array}\right)$$ 
be a $2\times 2$-minor of $B$, for some $3\leq i\neq j\leq n$. We remark that $B_{ij}$ cannot have a triangular shape. For example if $[e_2,e_i]=0$, then we could define a linear form $f\in \g^*$ such that $f$ is non-zero on $[e_1,e_i]$ and $[e_2,e_j]$ and zero outside $\CC[e_1,e_i]+\CC[e_2,e_j]$. Then $\{e_1,e_2,e_i\}$ would be linearly independent modulo $\g^f$, i.e. $\dim(\g/\g^f)\leq n-3$ contradicting the hypothesis that the index is $n-2$. 

If one of the rows of $B$ is zero, then $\af\oplus \CC e_i$ is an abelian ideal of codimension $1$ for $i=1$ or $i=2$, contradicting our assumption. Hence there exist $3\leq i,j\leq n$ such that $[e_1,e_i]\neq 0 \neq [e_2,e_j]$. We will show that after a suitable base change we can assume $i=j$.  

Assume $i\neq j$ and suppose first that $[e_1,e_i]$ and $[e_1,e_j]$ are linearly independent. 
Note that $[e_1,e_i]$ and $[e_2,e_i]$ are linearly independent, otherwise  if $[e_2,e_i]=\lambda [e_1,e_i]$ for some $\lambda \in \CC$, then after a base change $e_2 \leftarrow e_2-\lambda e_1$ we have $[e_2,e_i]=0$ and hence $B_{ij}$ has triangular shape which is not possible as mentioned before.
This allows us to define a linear form $f\in \g^*$ that is non-zero on $[e_1,e_j]$ and $[e_2,e_i]$, and zero on $[e_1,e_i]$. 
Then $\{e_1,e_2,e_i\}$ are linearly independent over $\g^f$ contradicting the assumption on the index. 
Hence we must have that $[e_1,e_j]$ and $[e_1,e_i]$ are linearly dependent, say $[e_1,e_j]=\lambda[e_1,e_i]$, $\lambda\in \CC$. After the base change $e_j \leftarrow e_j-\lambda e_i$ we have $[e_1,e_j]=0$, which gives the minor $B_{ij}$ a triangular shape since $[e_2,e_j]\neq 0$ and since $i\neq j$. But as said before, the minors $B_{ij}$ cannot have a triangular shape. Hence $i=j$ and without loss of generality we may assume $i=3$. Moreover we can change the basis of $\af$ such that $e_4=[e_1,e_3]$ and $e_5=[e_2,e_3]$. Note that since $\g/\af$ is abelian, $[e_1,e_2]\in \af$ and we are left with two cases:

\textbf{Case 1} If $[e_1,e_2] \not\in \langle e_3,e_4,e_5\rangle$, then change the basis of $\af$ such that $e_6=[e_1,e_2]$.
The Lie algebra $\g$ is then equal to the direct product $\g = \h_6 \times \af'$ where $\af'=\langle e_7, \ldots, e_n\rangle$. 

\textbf{Case 2} If $[e_1,e_2] \in \langle e_3,e_4,e_5\rangle$ then note first that $[e_1,e_2]\not\in\langle e_4,e_5\rangle$, since otherwise if 
$[e_1,e_2] = \alpha e_4 + \beta e_5$, for some $\alpha, \beta \in\CC$ we have after the base change $e_1\leftarrow e_1+\beta e_3$ and $e_2\leftarrow e_2-\alpha e_3$ that
$[e_1,e_2]=0$. Hence $\CC e_1 \oplus \CC e_2 \oplus \langle e_4,\ldots, e_n\rangle$ would be an abelian ideal of codimension $1$ - a contradiction to our hypothesis.
Thus $[e_1,e_2] = \alpha e_3 + \beta e_4 + \gamma e_5$ with $\alpha\neq 0$. After the basis change  
$$e_3 \leftarrow \alpha e_3 + \beta e_4 + \gamma e_5, \qquad e_4 \leftarrow \alpha^{-1} e_4, \qquad e_5 \leftarrow \alpha^{-1} e_5$$
we have $[e_1,e_2]=e_3$, $[e_1,e_3]=e_4$ and $[e_2,e_3]=e_5$. Hence $\g=\h_5 \times \af'$ where $\af' = \langle e_6, \ldots, e_n\rangle$ is abelian.
\end{proof}


\begin{proof}[Proof of the Main Theorem \ref{Main_Theorem}]
$(a)\Leftrightarrow (c)$ and$(b)\Leftrightarrow (c)$  follow from Proposition \ref{diamond-index}.
$(c)\Leftrightarrow (d)$ follows from Proposition \ref{listOfDiamond}.
\end{proof}

\section{Examples}
Finite dimensional Lie algebras $\g$ with an abelian ideal of codimension $1$ are in bijection with finite dimensional vector spaces $V$ and nilpotent endomorphisms $f:V\rightarrow V$. For such data one defines $\g=\CC e \oplus V$ and $[e,x]=f(x)$ for all $x\in V$. An example of this construction is given by the $n$-dimensional {\bf standard filliform} Lie algebra, which is the Lie algebra on the vector space $\mathcal{L}_n=\mathrm{span} \{e_1, \ldots, e_n\}$ such that $[e_1,e_i]=e_{i+1}$ for all $2\leq i < n$ and $[e_i,e_j]=0$ for $i,j\neq 1$.
Hence $\mathcal{L}_n$ provides an example of a non-abelian nilpotent Lie algebra $\g$ such that $U(\g)$ has property $(\diamond)$. The $3$-dimensional Heisenberg Lie algebra occurs as $\mathcal{L}_3$.

Given an even dimensional complex vector space $V=\CC^{2n}$ and an anti-symmetric bilinear form $\omega: V\times V \rightarrow \CC$, one defines the $2n+1$-dimensional {\bf Heisenberg Lie algebra} associated to $(V,\omega)$ as $\mathcal{H}_{2n+1}=V\oplus \CC h$ with $h$ being central and $[x,y]=\omega(x,y)h$ for all $x,y\in V$. Note that $\indx{\mathcal{H}_{2n+1}} = 1$. Thus $U(\mathcal{H}_{2n+1})$ satisfies $(\diamond)$ if and only if $n=1$, i.e. for $\mathcal{H}_{3}=\mathcal{L}_3$.

In \cite{HeisenbergLieSuperAlgebras} a finite dimensional Lie superalgebra $\g$ is called a {\bf Heisenberg Lie superalgebra} if it has a $1$-dimensional homogeneous center $\CC h= Z(\g)$ such that $[\g,\g] \subseteq Z(\g)$ and such that the associated homogeneous skew-supersymmetric bilinear form $\omega: \g\times \g \rightarrow \CC$ given by $[x,y]=\omega(x,y)h$ for all $x,y\in \g$ is non-degenerated when extended to $\g/Z(\g)$. On the other hand one can construct a Heisenberg Lie superalgebra on any finite-dimensional supersymplectic vector superspace $V$ with a homogeneous supersymplectic form $\omega$.

By \cite[page 73]{HeisenbergLieSuperAlgebras} if $\omega$ is even, i.e. $\omega(\g_0,\g_1)=0$, then $\g_0$ is a Heisenberg Lie algebra and if $\omega$ is odd, i.e. $\omega(\g_i,\g_i)=0$ for $i\in \{0,1\}$, then $\g_0$ is Abelian.
Hence $U(\g)$ satisfies $(\diamond)$ if and only if $\omega$ is odd or $\dim{\g_0}\leq 3$.

\bibliographystyle{amsplain}
\providecommand{\bysame}{\leavevmode\hbox to3em{\hrulefill}\thinspace}
\providecommand{\MR}{\relax\ifhmode\unskip\space\fi MR }
\providecommand{\MRhref}[2]{%
  \href{http://www.ams.org/mathscinet-getitem?mr=#1}{#2}
}
\providecommand{\href}[2]{#2}

\end{document}